\documentclass{amsart}

\usepackage[leqno]{amsmath}
\usepackage{amsthm}
\usepackage{latexsym,amsfonts,amssymb}
\usepackage[all]{xy} \SelectTips{eu}{}
\usepackage{hyperref}

\newcommand{\numberseries}{\bfseries}   
\newlength{\thmtopspace}                
\newlength{\thmbotspace}                
\newlength{\thmheadspace}               
\newlength{\thmindent}                  
\newtheoremstyle{bfupright head,slanted body}
                {\thmtopspace}{\thmbotspace}
                {\slshape}{\thmindent}{\bfseries}{.}{\thmheadspace}
                {{\numberseries \thmnumber{#2\;}}\thmnote{#3}}
\newtheoremstyle{fixed bf head,slanted body}
                {\thmtopspace}{\thmbotspace}{\slshape}
                {\thmindent}{\bfseries}{.}{\thmheadspace}
                {{\numberseries \thmnumber{#2\;}}\thmname{#1}\thmnote{ (#3)}}

\newtheoremstyle{fixed bf head,upright body}
                {\thmtopspace}{\thmbotspace}{\upshape}
                {\thmindent}{\bfseries}{.}{\thmheadspace}
                {{\numberseries \thmnumber{#2\;}}\thmname{#1}\thmnote{ (#3)}}
\newtheoremstyle{numbered paragraph}
                {\thmtopspace}{\thmbotspace}{\upshape}
                {\thmindent}{\upshape}{}{\thmheadspace}
                {{\numberseries \thmnumber{#2.}}}
\theoremstyle{bfupright head,slanted body}
\newtheorem{res}{}[section]             \newtheorem*{res*}{}
\theoremstyle{fixed bf head,slanted body}
\newtheorem{thm}[res]{Theorem}          \newtheorem*{thm*}{Theorem}
\newtheorem{prp}[res]{Proposition}      \newtheorem*{prp*}{Proposition}
\newtheorem{cor}[res]{Corollary}        \newtheorem*{cor*}{Corollary}
\newtheorem{lem}[res]{Lemma}            \newtheorem*{lem*}{Lemma}
\theoremstyle{fixed bf head,upright body}
\newtheorem{dfn}[res]{Definition}       \newtheorem*{dfn*}{Definition}
\newtheorem{rmk}[res]{Remark}           \newtheorem*{rmk*}{Remark}
          \newtheorem*{exa*}{Example}
\theoremstyle{numbered paragraph}
\newtheorem{ipg}[res]{}
\newlength{\thmlistleft}        
\newlength{\thmlistright}       
\newlength{\thmlistpartopsep}   
\newlength{\thmlisttopsep}      
\newlength{\thmlistparsep}      
\newlength{\thmlistitemsep}     
\newcounter{eqc}
\newenvironment{eqc}{\begin{list}{\upshape (\textit{\roman{eqc}})}%
    {\usecounter{eqc}%
      \setlength{\leftmargin}{\thmlistleft}%
      \setlength{\labelwidth}{\thmlistleft}%
      \setlength{\rightmargin}{\thmlistright}%
      \setlength{\partopsep}{\thmlistpartopsep}%
      \setlength{\topsep}{\thmlisttopsep}%
      \setlength{\parsep}{\thmlistparsep}%
      \setlength{\itemsep}{\thmlistitemsep}}}%
  {\end{list}}%
\newcommand{\eqclbl}[1]{{\upshape(\textit{#1})}}
\newcounter{prt}
\newenvironment{prt}{\begin{list}{\upshape (\alph{prt})}%
    {\usecounter{prt}%
      \setlength{\leftmargin}{\thmlistleft}%
      \setlength{\labelwidth}{\thmlistleft}%
      \setlength{\rightmargin}{\thmlistright}%
      \setlength{\partopsep}{\thmlistpartopsep}%
      \setlength{\topsep}{\thmlisttopsep}%
      \setlength{\parsep}{\thmlistparsep}%
      \setlength{\itemsep}{\thmlistitemsep}}}%
  {\end{list}}%

\newcounter{rqm}
  {\end{list}}%

  {\end{list}}%
\newenvironment{prf*}[1][Proof]{%
  \begin{proof}[\bf #1]
    \setcounter{equation}{0}
    }
  {\end{proof}
}
\newcommand{\proofofimp}[3][:]{\mbox{\eqclbl{#2}$\!\implies\!$\eqclbl{#3}#1}}

\newcommand{\pgref}[1]{\ref{#1}}

\renewcommand{\eqref}[1]{(\pgref{eq:#1})}

\numberwithin{equation}{res}

\newcommand{\fsp}[1][1.2]{f\hspace{#1pt}}

\newcommand{\setof}[3][\mspace{1mu}]{\{#1#2 \mid #3#1\}}
\newcommand{\kk}{\Bbbk}

\newcommand{\ZZ}{\mathbb{Z}}
\newcommand{\qtext}[1]{\quad\text{#1}\quad}

\newcommand{\qand}{\qtext{and}}

\newcommand{\deq}{\:=\:}
\newcommand{\dis}{\:\is\:}

\newcommand{\is}{\cong}
\newcommand{\qis}{\simeq}
\renewcommand{\le}{\leqslant}
\renewcommand{\ge}{\geqslant}
\newcommand{\lra}{\longrightarrow}
\newcommand{\xra}[2][]{\xrightarrow[#1]{\:#2\:}}
\newcommand{\qra}{\xra{\smash{\qis}}}
\newcommand{\Rop}{R^\circ}

\newcommand{\mapdef}[4][\rightarrow]{\nobreak{#2\colon #3 #1 #4}}

\newcommand{\Ker}[1]{\nobreak{\operatorname{Ker}#1}}
\newcommand{\Coker}[1]{\nobreak{\operatorname{Coker}#1}}
\renewcommand{\Im}[1]{\nobreak{\operatorname{Im}#1}}
\newcommand{\Cone}[1]{\nobreak{\operatorname{Cone}#1}}
\newcommand{\dif}[2][]{{\partial}^{#2}_{#1}}
\newcommand{\QQ}{\mathbb{Q}}
\newcommand{\Tsa}[2]{#2_{{\scriptscriptstyle\subset}#1}}
\newcommand{\Co}[2][]{\operatorname{C}_{#1}(#2)}
\renewcommand{\H}[2][]{\operatorname{H}_{#1}(#2)}
\newcommand{\Tha}[2]{#2_{{\scriptscriptstyle\le}#1}}
\newcommand{\Thb}[2]{#2_{{\scriptscriptstyle\ge}#1}}

\newcommand{\fd}[2][R]{\operatorname{fd}_{#1}#2}
\newcommand{\id}[2][R]{\operatorname{id}_{#1}#2}
\newcommand{\pd}[2][R]{\operatorname{pd}_{#1}#2}

\newcommand{\Hom}[3][R]{\operatorname{Hom}_{#1}(#2,#3)}

\newcommand{\Ext}[4][R]{\operatorname{Ext}_{#1}^{#2}(#3,#4)}
\newcommand{\tp}[3][R]{\nobreak{#2\otimes_{#1}#3}}

\newcommand{\Tor}[4][R]{\operatorname{Tor}^{#1}_{#2}(#3,#4)}
\newcommand{\Cat}[2]{{\mathsf{#2}}(#1)}
\newcommand{\splfR}[1][R]{\operatorname{splf}#1}

\newcommand{\D}[1][R]{\Cat{#1}{D}}

\newcommand{\del}{\partial}

\newcommand{\btp}[3][R]{\nobreak{#2\mathbin{\widebar{\otimes}}_{#1}#3}}

\newcommand{\ttp}[3][R]{\nobreak{#2\mathbin{\widetilde{\otimes}}_{#1}#3}}

\newcommand{\bTor}[4][R]{\widebar{\operatorname{Tor}}_{#2}^{#1}(#3,#4)}
\newcommand{\Ttor}[4][R]{\smash{\operatorname{\widehat{Tor}}}_%
  {#2}^{{#1}^{\phantom{|\mspace{-6mu}}}}(#3,#4)}
\newcommand{\Stor}[4][R]{\smash{\operatorname{\widetilde{Tor}}}_%
  {#2}^{{#1}^{\phantom{|\mspace{-6mu}}}}(#3,#4)}

\newcommand{\susp}{\operatorname{\mathsf{\Sigma}}}

\newcommand{\Gpd}[2][R]{\operatorname{Gpd}_{#1}#2}
\newcommand{\Gfd}[2][R]{\operatorname{Gfd}_{#1}#2}
\newcommand{\Gid}[2][R]{\operatorname{Gid}_{#1}#2}

\def\urltilda{\kern
  -.15em\lower .7ex\hbox{\~{}}\kern .04em} \makeatletter
\def\@nobreak@#1{\mathchoice%
  {\nobreakdef@\displaystyle\f@size{#1}}%
  {\nobreakdef@\nobreakstyle\tf@size{\firstchoice@false #1}}%
  {\nobreakdef@\nobreakstyle\sf@size{\firstchoice@false #1}}%
  {\nobreakdef@\nobreakstyle\ssf@size{\firstchoice@false #1}}%
  \check@mathfonts}%
\def\nobreakdef@#1#2#3{\hbox{{%
      \everymath{#1}%
      \let\f@size#2\selectfont%
      #3}}}%
\makeatother

\def\widebardisplay#1{%
  \setbox0=\hbox{$\displaystyle #1$} \dimen0=\wd0%
  \advance\dimen0 by -1pt
  \vbox{%
    \nointerlineskip%
    \moveright
    .25pt 
    \vbox{\hrule width \dimen0}%
    \nointerlineskip%
    \kern 2.25pt
    \box0%
  }%
}

\def\widebartext#1{%
  \setbox0=\hbox{$#1$} \dimen0=\wd0%
  \advance\dimen0 by -1pt
  \vbox{%
    \nointerlineskip%
    \moveright
    .25pt 
    \vbox{\hrule width \dimen0}%
    \nointerlineskip%
    \kern 1.75pt
    \box0%
  }%
}

\def\widebarscript#1{%
  \setbox0=\hbox{$\scriptstyle #1$} \dimen0=\wd0%
  \advance\dimen0 by -2pt
  \vbox{%
    \nointerlineskip%
    \moveright
    1pt 
    \vbox{\hrule width \dimen0}%
    \nointerlineskip%
    \kern .8pt
    \box0%
  }%
}

\def\widebarscriptscript#1{%
  \setbox0=\hbox{$\scriptscriptstyle #1$} \dimen0=\wd0%
  \advance\dimen0 by -2pt
  \vbox{%
    \nointerlineskip%
    \moveright
    1pt 
    \vbox{\hrule width \dimen0}%
    \nointerlineskip%
    \kern .6pt
    \box0%
  }%
}

\def\widebar#1{\mathchoice%
  {\widebardisplay{#1}}%
  {\widebartext{#1}}%
  {\widebarscript{#1}}%
  {\widebarscriptscript{#1}}%
}

\DeclareFontFamily{U}{mathx}{\hyphenchar\font45}
\DeclareFontShape{U}{mathx}{m}{n}{ <5> <6> <7> <8> <9> <10> <10.95>
  <12> <14.4> <17.28> <20.74> <24.88> mathx10 }{}
\DeclareSymbolFont{mathx}{U}{mathx}{m}{n}
\DeclareFontSubstitution{U}{mathx}{m}{n}
\DeclareMathAccent{\widecheck}{0}{mathx}{"71}
\DeclareMathAccent{\wideparen}{0}{mathx}{"75}

\def\bF{{\boldsymbol F}}
\def\bG{{\boldsymbol G}}
\def\bI{{\boldsymbol I}}
\def\bK{{\boldsymbol K}}

\def\bM{{\boldsymbol M}}
\def\bN{{\boldsymbol N}}
\def\bP{{\boldsymbol P}}
\def\bT{{\boldsymbol T}}
\def\bX{{\boldsymbol X}}

\def\bS{{\boldsymbol S}}

   \def\soft#1{\leavevmode\setbox0=\hbox{h}\dimen7=\ht0\advance
    \dimen7 by-1ex\relax\if t#1\relax\rlap{\raise.6\dimen7
    \hbox{\kern.3ex\char'47}}#1\relax\else\if T#1\relax
    \rlap{\raise.5\dimen7\hbox{\kern1.3ex\char'47}}#1\relax
    \else\if d#1\relax\rlap{\raise.5\dimen7\hbox{\kern.9ex
    \char'47}}#1\relax\else\if D#1\relax\rlap{\raise.5\dimen7
    \hbox{\kern1.4ex\char'47}}#1\relax\else\if l#1\relax
    \rlap{\raise.5\dimen7\hbox{\kern.4ex\char'47}}#1\relax
    \else\if L#1\relax\rlap{\raise.5\dimen7\hbox{\kern.7ex
    \char'47}}#1\relax\else\message{accent \string\soft
    \space #1 not defined!}#1\relax\fi\fi\fi\fi\fi\fi}

\begin{document}

\title{Homology theories and Gorenstein dimensions for complexes}

\author[L. Liang]{Li Liang}

\address{School of Mathematics and Physics, Lanzhou Jiaotong University, Lanzhou 730070, China}

\email{lliangnju@gmail.com}

\thanks{This research was partly supported by the National Natural Science Foundation of China (Grant Nos. 11761045 and 11971388), the Foundation of A Hundred Youth Talents Training Program of Lanzhou Jiaotong University, and the Natural Science Foundation of Gansu Province (Grant No. 18JR3RA113).}

\date{\today}



\begin{abstract}
In this paper, we first study the Gorenstein projective/flat dimension of complexes of modules. The relation between the Gorenstein projective/flat dimension for complexes and that for modules are investigated. Then we study Tate, stable and unbounded homology for complexes of modules. In the case of module arguments, we get some results that improve the known results.
\end{abstract}

\maketitle

\thispagestyle{empty}

\section*{Introduction}

\noindent
Tate (co)homology was initially defined for representations of finite groups. It was extended by Avramov and Martsinkovsky \cite{AM}  and  Veliche \cite{V} to modules/complexes of finite Gorenstein projective dimension. Christensen and Jorgensen further studied Tate homology for complexes of finite Gorenstein projective dimension in \cite{CJ}, where they proved a balanced result using pinched complexes: If $\bM$ is an $\Rop$-complex and $\bN$ is an $R$-complex, both of which are
bounded above and of finite Gorenstein projective dimension with $\bT \to \bP \to \bN$ a complete projective resolution of $\bN$, then for each $i\in\ZZ$ there is an isomorphism $\Ttor{i}{\bM}{\bN}\cong\H[i]{\tp{\bM}{\bT}}$. In Section \ref{6}, we reprove this result, and give another balanced result for Tate homology of complexes via flat objects. More precisely, we prove the next result; see Theorem \ref{cor:fb}.

\begin{res*}[Theorem A]
Let $\bM$ be a bounded above $\Rop$-complex of finite Gorenstein projective dimension and let $\bN$ be a bounded above $R$-complex of finite Gorenstein flat dimension with $(\bT, \bF)$ a Tate flat resolution. Then for each $i\in\ZZ$ there is an isomorphism
$$\Ttor{i}{\bM}{\bN}\cong\H[i]{\tp{\bM}{\bT}}.$$
\end{res*}

Here the definitions of Gorenstein projective dimension and Gorenstein flat dimension for complexes of modules were proposed by Veliche \cite{V} and Iacob \cite{I2}, respectively. In Section \ref{3}, we further study these Gorenstein homological dimensions. In particular, we prove the following result, which gives the relation between the Gorenstein projective/flat dimension for complexes and that for modules; see Theorems \ref{prp:relation} and \ref{prp:relation-gpd}.

\begin{res*}[Theorem B]
Let $\bM$ be a bounded $\Rop$-complex. Then the following hold:
  \begin{prt}
  \item There is an inequality
  $$\Gpd[\Rop]\bM \leq \max\{\Gpd[\Rop]M_i\ |\ i\in\ZZ\}+\sup \bM.$$
  In particular, if $\Gpd[\Rop]M_i<\infty$ for each $i\in\ZZ$ then $\Gpd[\Rop]\bM<\infty$.
  \item There is an inequality
  $$\Gfd[\Rop]\bM \leq \max\{\Gfd[\Rop]M_i\ |\ i\in\ZZ\}+\sup \bM.$$
  In particular, if $\Gfd[\Rop]M_i<\infty$ for each $i\in\ZZ$ then $\Gfd[\Rop]\bM<\infty$.
  \end{prt}
\end{res*}
We notice that the above inequalities can be strict: Let $P$ be a projective $\Rop$-module, and let $\bM$ be the acyclic complex with $M_0=M_{-1}=P$ and all other entries $0$. Then one has $\Gpd[\Rop]\bM=-\infty=\Gfd[\Rop]\bM$. However, the numbers $\max\{\Gpd[\Rop]M_i\ |\ i\in\ZZ\}+\sup \bM$ and $\max\{\Gfd[\Rop]M_i\ |\ i\in\ZZ\}+\sup \bM$ are zero.

\begin{equation*}
  \ast \ \ \ast \ \ \ast
\end{equation*}
\noindent
As a broad generalization of Tate homology to the realm of associative rings, stable homology was introduced by Vogel and Goichot \cite{G}, and further treated by Celikbas, Christensen, Liang and Piepmeyer \cite{CCLP1, CCLP2}, Emmanouil and Manousaki \cite{EM}, and Liang \cite{L2}. There are tight connections between stable homology functors $\widetilde{\mathrm{Tor}}$, absolute homology functors $\mathrm{Tor}$ and unbounded homology functors $\overline{\mathrm{Tor}}$; see \cite[2.5]{CCLP1}. Actually, there is an exact sequence of functors
$$\cdots \to \widetilde{\mathrm{Tor}}_{i} \to \mathrm{Tor}_{i} \to \overline{\mathrm{Tor}}_{i} \to \widetilde{\mathrm{Tor}}_{i-1} \to \cdots.$$

Sections \ref{4} and \ref{5} are devoted to studying stable and unbounded homology.

In Section \ref{4}, we give some developments related to Auslander's use of the transpose in terms of stable and unbounded homology modules. We prove the next result, which is contained in Theorem \ref{cor:sequence}.

\begin{res*}[Theorem C]
  Let $M$ be a Gorenstein flat $\Rop$-module. Then for each $R$-module $N$ there is an
  exact sequence
    \begin{equation*}
    0 \to \Stor{0}{M}{N} \to \tp{M}{N} \to \bTor{0}{M}{N}
    \to \Stor{-1}{M}{N} \to 0\;.
  \end{equation*}
\end{res*}

As an application, we prove in Corollary \ref{CJ} that if $R$ is a noetherian ring and $M$ is a finitely generated Gorenstein projective $\Rop$-module, then for each $R$-module $N$ there is an exact sequence
    \begin{equation*}
    0 \to \Ttor{0}{M}{N} \to \tp{M}{N} \to \Hom{\Hom[\Rop]{M}{R}}{N}
    \to \Ttor{-1}{M}{N} \to 0\;.
  \end{equation*}
This result was first proved by Christensen and Jorgensen over a commutative noetherian local ring;  see \cite[Proposition 6.3]{CJ}.

Finally, in Section \ref{5}, we further study unbounded homology for complexes of finite Gorenstein flat dimension. In the case of module arguments, we compare unbounded homology to relative homology.

\section{Preliminaries}\label{pre}

\noindent
We begin with some notation and terminology for use throughout this paper.

\begin{ipg}
Throughout this work, all rings are assumed to be associative algebras over a
commutative ring $\kk$. Let $R$ be a ring. By an $R$-module we mean a left $R$-module, and we refer to right $R$-modules as modules over the opposite ring $\Rop$. By an $R$-complex $\bM$ we mean a complex of $R$-modules as follows:
  \begin{equation*}
    \cdots \lra M_{i+1} \xra{\dif[i+1]{\bM}} M_i \xra{\dif[i]{\bM}}
    M_{i-1} \lra \cdots\: .
  \end{equation*}
We frequently (and without warning) identify $R$-modules with $R$-complexes concentrated in degree $0$. For an $R$-complex $\bM$, we set $\sup \bM=\sup\{i\in\ZZ\ |\ M_{i}\neq 0\}$ and $\inf \bM=\inf\{i\in\ZZ\ |\ M_{i}\neq 0\}$. An $R$-complex $\bM$ is \emph{bounded above} if $\sup \bM<\infty$, and it is \emph{bounded below} if $\inf \bM>-\infty$. An $R$-complex $\bM$ is \emph{bounded} if it is both bounded above and bounded below. For $n\in\mathbb{Z}$, the symbol $\susp^{n}\bM$ denotes the complex with $(\susp^{n}\bM)_{i}=M_{i-n}$ and $\del_{i}^{\susp^{n}\bM}=(-1)^{n}\del_{i-n}^{\bM}$ for all $i\in\ZZ$. We set $\susp \bM=\susp^{1}\bM$.

For $n\in\ZZ$, the symbol $\Co[n]{\bM}$ denotes the cokernel of
  $\del_{n+1}^{\bM}$, and $\H[n]{\bM}$ denotes the $n$th \emph{homology} of $\bM$, i.e., $\Ker\del_{n}^{\bM}/\Im\del_{n+1}^{\bM}$. An $R$-complex $\bM$ with $\H{\bM}=0$ is called \emph{acyclic,} and a morphism of $R$-complexes $\bM
\to \bN$ that induces an isomorphism $\H{\bM} \to \H{\bN}$ is called a \emph{quasi-isomorphism}. The symbol $\qis$ is used to decorate quasi-isomorphisms, and it is also used for isomorphisms in derived categories.

The symbol $\bM_{\leq n}$ denotes the subcomplex of $\bM$ with $(\bM_{\leq n})_i=M_i$ for $i\leq n$ and $(\bM_{\leq n})_i=0$ for $i > n$, and the symbol $\bM_{\ge n}$ denotes the quotient complex of $\bM$ with $(\bM_{\ge n})_i=M_i$ for $i\ge n$ and $(\bM_{\ge n})_i=0$ for $i < n$. The symbol $\Tsa{n}{\bM}$ denotes the quotient complex of $\bM$ with $(\Tsa{n}{\bM})_i=M_i$ for $i\leq n-1$, $(\Tsa{n}{\bM})_n=\Co[n]{\bM}$ and $(\Tsa{n}{\bM})_i=0$ for $i > n$.
\end{ipg}

\begin{ipg}
  \label{tensor}
  For an $\Rop$-complex $\bM$ and an $R$-complex $\bN$, the tensor product
  $\tp{\bM}{\bN}$ is the $\kk$-complex with the degree-$n$ term
  $(\bM\otimes_{R}\bN)_{n}=\coprod_{i\in\ZZ}(M_{i}\otimes_{R}N_{n-i})$ and the
  differential given by $\del(x\otimes y)=\dif[i]{\bM}(x)\otimes
  y+(-1)^{i}x\otimes\dif[n-i]{\bN}(y)$ for $x\in M_{i}$ and $y\in
  N_{n-i}$.
\end{ipg}

The next lemma is standard; we include a proof for the convenience of the reader.

\begin{lem}\label{lem:bounded}
Let $\bM$ be an $\Rop$-complex and $\bN$ an $R$-complex. The following assertions
  hold for each $n\in\ZZ$.
    \begin{prt}
  \item If\, $\sup \bN\leq s$, then $\H[n]{\tp{\bM}{\bN}}=\H[n]{\tp{\bM_{\geq n-s-1}}{\bN}}$.
  \item If\, $\inf \bN\geq t$, then $\H[n]{\tp{\bM}{\bN}}=\H[n]{\tp{\bM_{\leq n-t+1}}{\bN}}$.
  \end{prt}
\end{lem}
\begin{prf*}
(a) For an integer $k\geq n-1$, one has
  \begin{align*}
    (\tp{\bM}{\bN})_{k}
    & = \coprod_{i\in\ZZ}(M_{i}\otimes_{R}N_{k-i})\\
    & = \coprod_{i\geq k-s}(M_{i}\otimes_{R}N_{k-i})\\
    & = \coprod_{i\geq k-s}((\bM_{\geq n-s-1})_{i}\otimes_{R}N_{k-i})\\
    & = \coprod_{i\in\ZZ}((\bM_{\geq n-s-1})_{i}\otimes_{R}N_{k-i})\\
    & = (\tp{\bM_{\geq n-s-1}}{\bN})_{k}\:,
  \end{align*}
where the second and fourth equalities hold as $\sup \bN\leq s$. Hence the equality $\H[n]{\tp{\bM}{\bN}}=\H[n]{\tp{\bM_{\geq n-s-1}}{\bN}}$ follows.

The proof of part (b) is similar.
\end{prf*}

\begin{ipg}
  \label{res}
  A complex $\bP$ of projective $R$-modules is called \emph{semi-projective} if the functor $\Hom{\bP}{-}$ preserves
  quasi-isomorphisms. A complex $\bI$ of injective $R$-modules is called
  \emph{semi-injective} if the
  functor $\Hom{-}{\bI}$ preserves quasi-isomorphisms. A complex
  $\bF$ of flat $R$-modules is called \emph{semi-flat} if the
  functor $\tp{-}{\bF}$ preserves quasi-isomorphisms.

  A \emph{semi-projective resolution} of $\bM$ is a quasi-isomorphism
    $\mapdef{\pi}{\bP}{\bM}$, where $\bP$ is a semi-projective complex. Every $R$-complex $\bM$ has a semi-projective resolution. Moreover, $\pi$ can be chosen surjective, and if $\inf\H{\bM}>-\infty$ then $\bP$ can be chosen such that $\inf \bP=\inf\H{\bM}$; see Avramov and Foxby \cite{AF}. Dually,
  A \emph{semi-injective resolution} of $\bM$ is a quasi-isomorphism
    $\mapdef{\iota}{\bM}{\bI}$, where $\bI$ is a semi-injective complex. Every $R$-complex $\bM$ has a semi-injective resolution. Moreover, $\iota$ can be chosen injective, and if $\sup\H{\bM}<\infty$ then $\bI$ can be chosen such that $\sup \bI=\sup\H{\bM}$. A \emph{semi-flat replacement} of $\bM$ is an isomorphism $\bF
  \qis \bM$ in the derived category, where $\bF$ is
  a semi-flat complex. Every complex has a
  semi-projective resolution and hence a semi-flat replacement.
\end{ipg}

\begin{ipg}
Let $\bM$ be an $\Rop$-complex with $\bP\qra \bM$ a semi-projective resolution, and let $\bN$ be an $R$-complex. For all $i\in\ZZ$ the modules
  \begin{equation*}
    \Tor{i}{\bM}{\bN} = \H[i]{\tp{\bP}{\bN}}
  \end{equation*}
make up the \emph{absolute homology} of $\bM$ and $\bN$ over $R$. It is clear that the definition of $\Tor{i}{\bM}{\bN}$ is functorial and homological in either argument.
\end{ipg}

\begin{ipg}
  \label{ipg:gproj}
  An acyclic complex $\bT$ of projective $\Rop$-modules is called
  \emph{totally acyclic} if $\Hom[\Rop]{\bT}{P}$ is acyclic for every
  projective $\Rop$-module $P$. An $\Rop$-module $G$ is called \emph{Gorenstein projective} if there exists a totally acyclic complex $\bT$ of projective $\Rop$-modules such that $G\is\Coker(T_1\to T_0)$.  A
  \emph{complete projective resolution} of an $\Rop$-complex $\bM$ is a
  diagram $\bT \xra{\tau} \bP \qra \bM$, where $\bT$ is a totally acyclic
  complex of projective $\Rop$-modules, $\bP\qra \bM$ is a semi-projective
  resolution, and $\tau_i$ is an isomorphism for $i \gg 0$. The \emph{Gorenstein projective dimension} of $\bM$ is defined by
    \begin{align*}
    \Gpd[\Rop]{\bM}
    & \deq \inf\left\{ g \in \ZZ \:
      \left|
        \begin{array}{c}
          \bT \xra{\tau} \bP \qra \bM\\
          \text{is a complete projective resolution with}\\
          \text{$\tau_i: T_i \to P_i$ isomorphism for each $i\geq g$}
        \end{array}
      \right.
    \right\};
  \end{align*}
  see \cite[Definition 3.1]{V}. It is clear that $\bM$ has finite Gorenstein projective dimension if and only if $\bM$ has a complete projective resolution.

  Let $\bM$ be an $\Rop$-complex that has a complete projective resolution $\bT
  \to \bP \to \bM$. From \cite[(2.4)]{CJ}, for an $R$-complex $\bN$ and $i\in\ZZ$, the $i$th \emph{Tate homology} of $\bM$
  and $\bN$ over $R$ is defined as
  $$\Ttor{i}{\bM}{\bN}=\H[i]{\tp{\bT}{\bN}}.$$
\end{ipg}

\section{Gorenstein homological dimensions of complexes}\label{3}

\noindent
In this section, we study the Gorenstein flat/projective dimension for complexes of modules. We start by recalling the following definitions.

\begin{ipg}
  \label{ipg:gproj}
  An acyclic complex $\bT$ of flat $\Rop$-modules is called
  \emph{F-totally acyclic} if $\tp{\bT}{E}$ is acyclic for every
  injective $R$-module $E$. An $\Rop$-module $L$ is called \emph{Gorenstein flat} if there exists a F-totally acyclic complex $\bT$ of flat $\Rop$-modules such that $L\is\Coker(T_1\to T_0)$. A \emph{Tate flat resolution} of an $\Rop$-complex $\bM$ is a pair $(\bT,\bF)$
  where $\bT$ is an F-totally acyclic complex of flat $\Rop$-modules and $\bF\simeq
  \bM$ is a semi-flat replacement with $\Thb{g}{\bT} \is \Thb{g}{\bF}$ for some
  $g\in\mathbb{Z}$; see Liang \cite{L1}. If furthermore, there exists a morphism $\tau: \bT \to \bF$
  such that $\tau_{i}$ is an isomorphism for each $i\geq g$,
  then the Tate flat resolution $(\bT,\bF)$ is said to be a \emph{complete flat resolution} of $\bM$.
\end{ipg}

The next definition can be found in Christensen, K\"{o}ksal and Liang \cite{CKL}.

\begin{dfn}\label{dfn:gfd}
  Let $\bM$ be an $\Rop$-complex. The \emph{Gorenstein flat dimension} of
  $\bM$ is given by
  \begin{equation*}
    \Gfd[\Rop]{\bM} \deq \inf\setof{g\in\ZZ}{(\bT,\bF)\
      \text{is a Tate flat resolution of $\bM$ with $\Thb{g}{\bT}
        \is \Thb{g}{\bF}$}}\:.
  \end{equation*}
\end{dfn}

\begin{rmk}\label{rmk:gfd}
From a recent result by \v{S}aroch and \v{S}{\soft{t}}ov{\'{\i}}{\v{c}}ek \cite[Theorem 4.11]{SS}, all rings are GF-closed in the terminology of Bennis \cite{Be}, and so Definition \ref{dfn:gfd} agrees with Iacob's definition \cite[Definition 3.2]{I2}; see \cite[Remark 5.13]{CKL}. Thus from \cite[Remark 3]{I2}, it extends the definitions in Holm \cite[3.9]{H1} and Christensen, Frankild and Holm \cite[1.9]{CFH} of Gorenstein flat dimension for modules and complexes with bounded below homology.
\end{rmk}

\begin{rmk}\label{rmk:gfd-acyclic}
  Let $\bM$ be an $\Rop$-complex. It is clear that $\Gfd[\Rop]{\bM}<\infty$ if and only if $\bM$ admits a Tate flat resolution, and $\Gfd[\Rop]{\bM}=-\infty$ if and only if $\bM$ is acyclic.
\end{rmk}

The next result provides complete flat resolutions, which is for use in Section \ref{5}.

\begin{thm}\label{prp:resolution}
Let $R$ be a left coherent ring and $g\in\ZZ$. For an $\Rop$-complex $\bM$ the following conditions
  are equivalent.
  \begin{eqc}
  \item $\Gfd[\Rop]{\bM}\leq g$.
  \item $\bM$ admits a complete flat resolution $\xymatrix@C=0.5cm{
    \bT \ar[r]^{\tau} & \bF}$ such that $\tau_{i}$ is split surjective for each $i\in\ZZ$ and $\tau_{i}=\id[F_{i}]$ for all $i\geq g$.
  \end{eqc}
\end{thm}
\begin{prf*}
The implication \proofofimp[]{ii}{i} is clear.

 \proofofimp{i}{ii} Let $\Gfd[\Rop]{\bM}\leq g$. By \cite[Proposition 5.7]{CKL} there exists a semi-flat replacement $\bF\simeq \bM$ such that $\Co[i]\bF$ is Gorenstein flat for each $i\geq g$, and $\H[i]{\bF}=0$ for all $i>g$. Since $R$ is left coherent, there is a flat preenvelope $f: \Co[g]\bF\to G_{g-1}$; see Xu \cite[Theorem 2.5.1]{X}. Note that $\Co[g]\bF$ is Gorenstein flat, then $f$ is a monomorphism, and so there is an exact sequence
   \begin{equation}
    \label{sequence1}
    \xymatrix@C=0.5cm{0 \ar[r] & \Co[g]\bF \ar[r]^{f} & G_{g-1} \ar[r] & C_{g-1} \ar[r] & 0}
  \end{equation}
 of $\Rop$-modules, which is $\Hom[\Rop]{-}{F}$ exact for each flat $\Rop$-module $F$. For each injective $R$-module $I$, one has $$\Hom[\ZZ]{\Tor{1}{C_{g-1}}{I}}{\QQ/\ZZ}\is \Ext[\Rop]{1}{C_{g-1}}{\Hom[\ZZ]{I}{\QQ/\ZZ}}=0,$$
 where the equality holds since $\Hom[\ZZ]{I}{\QQ/\ZZ}$ is a flat and cotorsion $\Rop$-module, and $f$ is a flat preenvelope. Hence $\Tor{1}{C_{g-1}}{I}=0$, and so $C_{g-1}$ is Gorenstein flat by \cite[Proposition 3.8]{H1}. Thus the sequence (\ref{sequence1}) is $\tp{-}{I}$ exact for each injective $R$-module $I$.

 Continue the above process, one gets an exact sequence
    \begin{equation}
    \label{sequence2}
    0 \to \Co[g]\bF \to G_{g-1} \to G_{g-2} \to \cdots
    \end{equation}
 of $\Rop$-modules with each $G_{i}$ flat, such that the sequence is $\tp{-}{I}$ and $\Hom[\Rop]{-}{F}$ exact for each injective $R$-module $I$ and flat $\Rop$-module $F$. So we have the following commutative diagram:
   \begin{equation*}
    \xymatrix{
      0 \ar[r]
      & \Co[g]\bF \ar[r] \ar@{=}[d]
      & G_{g-1} \ar[r] \ar[d]
      & G_{g-2} \ar[r] \ar[d]
      & \cdots\\
      0 \ar[r]
      & \Co[g]\bF \ar[r]
      & F_{g-1} \ar[r]
      & F_{g-2} \ar[r]
      & \cdots\:.
    }
  \end{equation*}
  Since $\Co[i]\bF$ is Gorenstein flat for each $i\geq g$, and $\H[i]{\bF}=0$ for all $i>g$, the sequence
  \begin{equation}
  \label{sequence3}
  \cdots \to F_{g+1} \to F_{g} \to \Co[g]{\bF} \to 0
  \end{equation}
  is exact, and it is $\tp{-}{I}$ exact for each injective $R$-module $I$.

  Assembling the sequences (\ref{sequence2}) and (\ref{sequence3}), one gets an F-totally acyclic complex $\bT'$ of flat $\Rop$-modules and a morphism $\alpha': \bT'\to \bF$ such that $\alpha'_{i}=\id[F_{i}]$ for each $i\geq g$.

  Set $\bT''=\susp^{-1}\Cone(\id[\bF_{<g}])$. Then $\bT''$ is a contractible complex and there is a degree-wise split surjective morphism $\tau: \bT'' \to \bF_{<g}$. Let $\alpha''=\varepsilon\tau: \bT'' \to \bF$, where $\varepsilon: \bF_{<g} \to \bF$ is the natural morphism. Then $\alpha''_{i}$ is split surjective for each $i<g$ and $\alpha''_{i}=0$ for each $i\geq g$. Let $\bT=\bT'\oplus \bT''$ and $\tau=(\alpha', \alpha''): \bT \to \bF$. Then $\bT$ is an F-totally acyclic complex of flat $\Rop$-modules, and $\tau_{i}$ is split surjective for each $i<g$ and $\tau_{i}=\id[F_{i}]$ for each $i\geq g$.
\end{prf*}

The following result complements the work of Iacob \cite{I2}, which is an analogue of \cite[Theorem 3.8]{V} for Gorenstein flat dimension.

\begin{prp}\label{prp:gfd}
Let $\bM$ be an $\Rop$-complex of finite Gorenstein flat dimension. Then one has
  \begin{align*}
    \Gfd[\Rop]{\bM}
    & \deq \sup\left\{ n \in \ZZ \:
      \left|
        \begin{array}{c}
          \Tor{n+\sup\H{\bN}}{\bM}{\bN}\neq 0\\
          \text{for some $R$-complex $\bN$ with}\\
          \text{$\id{\bN}<\infty$ and $\sup \H{\bN}<\infty$}
        \end{array}
      \right.
    \right\}\\
    & \deq \sup\left\{ n \in \ZZ \:
      \left|
        \begin{array}{c}
          \Tor{n}{\bM}{E}\neq 0\\
          \text{for some injective $R$-module $E$}
        \end{array}
      \right.
    \right\}
  \end{align*}
\end{prp}
\begin{prf*}
Let $\Gfd[\Rop]{\bM}=g<\infty$. If $g=-\infty$, then $M$ is acyclic; see Remark \ref{rmk:gfd-acyclic}. Thus $\Tor{n}{\bM}{-}=0$ for each $n\in\ZZ$, and so one gets the equalities in the statement. We let $g\in \ZZ$, and let $s$ (resp., $s'$) denote the number on the right side of the first (resp., second) equality in the statement. Then there exists a Tate flat resolution $(\bT, \bF)$ such that $\bT$ is an F-totally acyclic complex of flat $\Rop$-modules, $\bF\simeq \bM$ is a semi-flat replacement and $\bT_{\geq g} \is \bF_{\geq g}$. For each $R$-complex $\bN$ with $\id{\bN}<\infty$ and $\sup \H{\bN}<\infty$, choose a semi-injective resolution $\bN \qra \bI$ such that $\sup \bI =\sup\H{\bN}$ and $\inf \bI>-\infty$; see \cite[Theorem 2.4.I]{AF}. Then for all $i>g+\sup\H{\bN}$, one has
  \begin{align*}
    \Tor{i}{\bM}{\bN}
    & \dis \H[i]{\tp{\bF}{\bN}}\\
    & \dis \H[i]{\tp{\bF}{\bI}}\\
    & \deq \H[i]{\tp{\bF_{\geq i-\sup\H{N}-1}}{\bI}}\\
    & \dis \H[i]{\tp{\bT_{\geq i-\sup\H{N}-1}}{\bI}}\\
    & \deq \H[i]{\tp{\bT}{\bI}}\\
    & \deq 0\:,
  \end{align*}
where the first isomorphism follows from Christensen, Foxby and Holm \cite[Corollary 7.4.15]{CFoH}, the first and second equalities hold by Lemma \ref{lem:bounded}, and the last equality follows from \cite[Lemma 2.13]{CFH}. Thus one gets $g\geq s$. Obviously, $s\geq s'$, so it remains to prove $s'\geq g$.

Since $\Gfd[\Rop]{\bM}<\infty$, one has $\Gid{\Hom[\ZZ]{\bM}{\QQ/\ZZ}}<\infty$; see \cite[Proposition 5.8]{CKL}.
On the other hand, for each injective $R$-module $E$ and all $i>s'$, one has
  \begin{align*}
    \Ext{i}{E}{\Hom[\ZZ]{\bM}{\QQ/\ZZ}}
    & \dis \Hom[\ZZ]{\Tor{i}{\bM}{E}}{\QQ/\ZZ}\\
    & \deq 0\:,
  \end{align*}
so $\Gid{\Hom[\ZZ]{\bM}{\QQ/\ZZ}}\leq s'$ by Asadollahi and Salarian \cite[Theorem 2.4]{AS}. Thus by \cite[Theorem 2.3]{AS}, for all $i>s'$ one has
$$\Hom[\ZZ]{\H[i]{\bM}}{\QQ/\ZZ} \is \H[-i]{\Hom[\ZZ]{\bM}{\QQ/\ZZ}}=0,$$
and hence $\H[i]{\bM}=0$. That is, $\sup\H{\bM}\leq s'$. Since $\susp^{-s'}\bF_{\geq s'} \qra \Co[s']{\bF}$ is a flat resolution, for each injective $R$-module $E$ and all $i\geq 1$ one has
  \begin{align*}
    \Tor{i}{\Co[s']{\bF}}{E}
    & \dis \H[i]{\tp{\susp^{-s'}\bF_{\geq s'}}{E}}\\
    & \dis \H[i+s']{\tp{\bF_{\geq s'}}{E}}\\
    & \deq \H[i+s']{\tp{\bF}{E}}\\
    & \dis \Tor{i+s'}{\bM}{E}\\
    & \deq 0\:,
  \end{align*}
where the first equality follows from Lemma \ref{lem:bounded}, as $\sup E\leq i-1$. We notice that $\Gfd[\Rop]{\Co[s']{\bF}}<\infty$ and all rings are GF-closed; see \cite[Theorem 4.11]{SS}. Then $\Co[s']{\bF}$ is Gorenstein flat by \cite[Theorem 2.8]{Be}, and so one has $s'\geq g$ by \cite[Proposition 5.7]{CKL}.
\end{prf*}

\begin{cor}
Let $0 \to \bM' \to \bM \to \bM'' \to 0$ be an exact sequence of $\Rop$-complexes. Then there is an inequality $$\Gfd[\Rop]{\bM}\leq \max\{\Gfd[\Rop]{\bM'}, \Gfd[\Rop]{\bM''}\}.$$
\end{cor}
\begin{prf*}
We may assume that $\max\{\Gfd[\Rop]{\bM'}, \Gfd[\Rop]{\bM''}\}<\infty$. Since all rings are GF-closed by \cite[Theorem 4.11]{SS}, one has $\Gfd[\Rop]{\bM}<\infty$ by \cite[Proposition 3.4]{I2}. For each injective $R$-module $E$ there is an exact sequence
$$\cdots \to \Tor{n}{\bM'}{E} \to \Tor{n}{\bM}{E} \to \Tor{n}{\bM''}{E} \to \cdots.$$
Thus one gets $\Gfd[\Rop]{\bM}\leq \max\{\Gfd[\Rop]{\bM'}, \Gfd[\Rop]{\bM''}\}$ by Proposition \ref{prp:gfd}.
\end{prf*}

In the following we investigate the relation between the Gorenstein flat/projective dimension
for complexes of modules and that for modules.

\begin{lem}\label{lem:gfd}
Let $\bM$ be a bounded $\Rop$-complex such that $\Gfd[\Rop]M_i<\infty$ for each $i\in\ZZ$. Then there exists an exact sequence $0 \to \bK \to \bG \to \bM \to 0$ of $\Rop$-complex such that the following conditions hold:
  \begin{prt}
  \item $\bG$ is a bounded complex of Gorenstein flat $\Rop$-modules with $G_i$ flat for $i>\sup \bM$, $\inf \bG=\inf \bM$ and $\sup \bG \leq \max\{\Gfd M_i\ |\ i\in\ZZ\}+\sup \bM$.
  \item $\bK$ is a bounded acyclic $\Rop$-complex with each $\Co[i]\bK$ cotorsion $\Rop$-module of finite flat dimension.
  \end{prt}
\end{lem}
\begin{prf*}
Set $g=\max\{\Gfd[\Rop]M_i\ |\ i\in\ZZ\}$. One has $g<\infty$ since $\bM$ is a bounded $\Rop$-complex.  Without loss of generality, we may assume that $\inf \bM=0$ and $\sup \bM=s$. We argue by induction on $s$. If $s=0$, then one has $\bM=M_0$ with $\Gfd[\Rop]M_0 =g$. It follows from \cite[Corollary 4.12]{SS} that the class $\mathcal{GF}$ of Gorenstein flat $\Rop$-modules and the class $\mathcal{FC}$ of flat and cotorsion $\Rop$-modules are closed under extensions and direct summands, and $\mathcal{FC}$ is a cogenerator for $\mathcal{GF}$ in the sense of Auslander and Buchweitz \cite{AB}. Thus there is an exact sequence
$$\xymatrix@C=0.5cm{0 \ar[r] & K_1 \ar[r] & G_0 \ar[r]^{\delta_0} & M_0 \ar[r] & 0}$$
such that $G_0$ is a Gorenstein flat and $K_1$ is cotorsion with $\fd[\Rop]K_1 = g-1$ by \cite[Theorem 1.1]{AB}. For the module $K_1$, it follows from Bican, El Bashir and Enochs \cite[Proposition 2]{BEE} that there is an exact sequence
$$\xymatrix@C=0.5cm{0 \ar[r] & K_2 \ar[r] & F_1 \ar[r]^{\delta_1} & K_1 \ar[r] & 0}$$
such that $F_1$ is flat and $K_2$ is cotorsion with $\fd[\Rop]K_2 = g-2$. Continue this process, one gets an exact sequence
$$\xymatrix@C=0.5cm{0 \ar[r] & F_g \ar[r] & F_{g-1} \ar[r]^{\delta_{g-1}} & \cdots \ar[r] & F_1 \ar[r]^{\delta_1} & G_0 \ar[r]^{\delta_0} & M_0 \ar[r] & 0}$$
such that $G_0$ is Gorenstein flat, each $F_i$ is flat and each $\Ker\delta_i$ is cotorsion with $\fd[\Rop]\Ker\delta_i<\infty$. Let
$$\bG=\cdots \to 0 \to F_g \to \cdots \to F_1 \to G_0 \to 0 \to \cdots.$$
Then $\bG$ satisfies the condition (a) in the statement, and there is a surjective morphism $\alpha: \bG\to \bM$. Let $\bK=\Ker\alpha$. It is clear that $\bK$ satisfies the condition (b) in the statement.

We let $s>0$. Then there exists a morphism $f: \susp^{s-1}M_s \to \Tha{s-1}{\bM}$. By induction hypothesis there exists exact sequences $0 \to \bK' \to \bG' \to \susp^{s-1}M_s \to 0$ and $0 \to \bK'' \to \bG'' \to \Tha{s-1}{\bM} \to 0$ such that the following conditions hold:
  \begin{prt}
  \item $\bK'$ and $\bK''$ are bounded acyclic $\Rop$-complexes with each $\Co[i]{\bK'}$ and $\Co[i]{\bK''}$ cotorsion $\Rop$-modules of finite flat dimension.
  \item $\bG'$ and $\bG''$ are bounded complexes of Gorenstein flat $\Rop$-modules with $\inf \bG'=s-1$, $\inf \bG''=0$, $\sup \bG'\leq g+s-1$ and $\sup \bG''\leq g+s-1$.
  \item $G'_i$ and $G''_i$ are flat for $i>s-1$.
  \end{prt}
From Liang \cite[Lemma 2.7(2)]{L20}, $\bG'$ is a Gorenstein flat $\Rop$-complex in the sense of Enochs and Garc\'{i}a Rozas \cite{EG}. So it follows from \cite[Lemma 2.14(1)]{L20} that $\Ext[\mathcal{C}]{1}{\bG'}{\bK''}=0$ \footnote{The symbol $\Hom[\mathcal{C}]{\bM}{\bN}$ denotes the set of morphisms of $R$-complexes from $\bM$ to $\bN$, and $\Ext[\mathcal{C}]{i}{-}{-}$ is the right derived functor of $\Hom[\mathcal{C}]{-}{-}$.}. Thus there is a commutative diagram of $\Rop$-complexes
$$\xymatrix{
  0 \ar[r] & \bK' \ar[d]_{h} \ar[r]^{\alpha'} & \bG' \ar[d]_{g} \ar[r]^{\beta'\ \ \ } & \susp^{s-1}M_s \ar[d]_{f} \ar[r] & 0 \\
  0 \ar[r] & \bK'' \ar[r]^{\alpha''} & \bG'' \ar[r]^{\beta''\ \ \ } & \Tha{s-1}{\bM} \ar[r] & 0 .}$$
This yields an exact sequence
$$\xymatrix{0 \ar[r] & \Cone h \ar[r]^{\alpha} & \Cone g \ar[r]^{\beta} & \Cone f \ar[r] & 0}$$
of $\Rop$-complexes, where $\alpha_i=\left(\begin{matrix} \alpha''_i & 0 \\ 0 & \alpha'_{i-1} \end{matrix}\right)$ and $\beta_i=\left(\begin{matrix} \beta''_i & 0 \\ 0 & \beta'_{i-1} \end{matrix}\right)$ for each $i\in\ZZ$. We observe that $\Cone f=\bM$. Let $\bG=\Cone g$ and $\bK=\Cone h$. Then $\bG$ and $\bK$ satisfy the conditions (a) and (b) in the statement, respectively.
\end{prf*}

\begin{thm}\label{prp:relation}
Let $\bM$ be a bounded $\Rop$-complex. Then there is an inequality
$$\Gfd[\Rop]\bM \leq \max\{\Gfd[\Rop]M_i\ |\ i\in\ZZ\}+\sup \bM.$$
In particular, if $\Gfd[\Rop]M_i<\infty$ for each $i\in\ZZ$ then $\Gfd[\Rop]\bM<\infty$.
\end{thm}
\begin{prf*}
One may assume that $\bM\neq 0$. Let $\sup \bM=s\in\ZZ$. We assume that $\max\{\Gfd[\Rop]M_i\ |\ i\in\ZZ\}=g<\infty$. By Lemma \ref{lem:gfd}, there exists an exact sequence $\xymatrix@C=0.5cm{0 \ar[r] & \bK \ar[r] & \bG \ar[r]^{\pi\ } & \bM \ar[r] & 0 }$ such that $\bK$ is a bounded acyclic $\Rop$-complex with each $\Co[i]\bK$ cotorsion $\Rop$-module of finite flat dimension, and $\bG$ is a bounded complex of Gorenstein flat $\Rop$-modules with $\sup \bG \leq g+s$. Fix a semi-projective resolution $\pi': \bP \qra \bM$ such that $\bP$ is bounded below. Since $\bP$ is a semi-flat $\Rop$-complex, one has $\Ext[\mathcal{C}]{1}{\bP}{\bK}=0$ by Gillespie \cite[Proposition 3.6]{Gi}, and so there is a morphism $\alpha: \bP \to \bG$ such that $\pi\alpha=\pi'$. We observe that $\pi$ and $\pi'$ are quasi-isomorphisms. Then $\alpha$ is a quasi-isomorphism. Let $\bX=\susp^{-1}\Cone\alpha$. Then $\bX$ is an bonded below acyclic complex of Gorenstein flat $\Rop$-modules with $X_i=P_i$ for $i\geq g+s$, and so $\H[i]{\bP}=\H[i]{\bX}=0$ for $i>g+s$, and $\Co[g+s]{\bP}=\Co[g+s]{\bX}$ is Gorenstein flat by \cite[Theorem 2.3]{Be}. Thus one has $\Gfd \bM\leq g+s$ by \cite[Proposition 5.12]{CKL}.
\end{prf*}

\begin{ipg}\label{ipg:bgfd}
For bounded complexes, the notion of Gorenstein flat dimension in this paper is compatible with the one given in \cite[1.9]{CFH}; see Remark \ref{rmk:gfd}. That is, if $\bM$ is a bounded $\Rop$-complex then one has
  \begin{align*}
    \Gfd[\Rop]{\bM}
    & \deq \inf\left\{\sup\bG \:
      \left|
        \begin{array}{c}
          \bG\qis\bM\ \text{is a bounded below complex}\\
          \text{of Gorenstein flat $\Rop$-modules}
        \end{array}
      \right.
    \right\}.
  \end{align*}
\end{ipg}

Next we give an alternate proof of Theorem \ref{prp:relation}.

\begin{prf*}[Alternate proof of Theorem \ref{prp:relation}]
One may assume that $\bM\neq 0$. Let $\sup \bM=s\in\ZZ$. We assume that $\max\{\Gfd[\Rop]M_i\ |\ i\in\ZZ\}=g<\infty$. By Lemma \ref{lem:gfd}, there exists an exact sequence $\xymatrix@C=0.5cm{0 \ar[r] & \bK \ar[r] & \bG \ar[r]^{\pi\ } & \bM \ar[r] & 0 }$ such that $\bK$ is an acyclic $\Rop$-complex and $\bG$ is a bounded complex of Gorenstein flat $\Rop$-modules with $\sup \bG \leq g+s$. Thus $\pi: \bG\to \bM$ is a quasi-isomorphism, and so by the equality in \ref{ipg:bgfd} one has $\Gfd[\Rop]{\bM}\leq\sup \bG \leq g+s$.
\end{prf*}

We end this section with the following two results, which are proved like Lemma \ref{lem:gfd} and Theorem \ref{prp:relation}.

\begin{lem}\label{lem:gpd}
Let $\bM$ be a bounded $\Rop$-complex such that $\Gpd[\Rop]M_i<\infty$ for each $i\in\ZZ$. Then there exists an exact sequence $0 \to \bK \to \bG \to \bM \to 0$ of $\Rop$-complex such that the following conditions hold:
  \begin{prt}
  \item $\bG$ is a bounded complex of Gorenstein projective $\Rop$-modules with $G_i$ projective for $i>\sup \bM$, $\sup \bG \leq \max\{\Gpd M_i\ |\ i\in\ZZ\}+\sup \bM$ and $\inf \bG=\inf \bM$.
  \item $\bK$ is a bounded acyclic $\Rop$-complex such that each $\Co[i]\bK$ has finite projective dimension.
  \end{prt}
\end{lem}

\begin{thm}\label{prp:relation-gpd}
Let $\bM$ be a bounded $\Rop$-complex. Then there is an inequality
$$\Gpd[\Rop]\bM \leq \max\{\Gpd[\Rop]M_i\ |\ i\in\ZZ\}+\sup \bM.$$
In particular, if $\Gpd[\Rop]M_i<\infty$ for each $i\in\ZZ$ then $\Gpd[\Rop]\bM<\infty$.
\end{thm}

\section{A balanced result for Tate homology of complexes}\label{6}
\noindent
In this section, we give a balanced result for Tate homology of complexes of modules.
\begin{lem}\label{lem:syzygy}
Let $\bT$ be an acyclic complex of flat $\Rop$-modules, and let $\bN$ be a bounded above $R$-complex with $\bF'\qis\bN$ a semi-flat replacement. Then for each $n\geq\sup \bN$ there is an isomorphism $\tp{\bT}{\bN}\simeq \susp^{n}(\tp{\bT}{\Co[n]{\bF'}})$ in $\D[\kk]$.
\end{lem}
\begin{prf*}
Since $n\geq\sup\H{\bF'}$, there is a quasi-isomorphism $\xymatrix@C=0.5cm{\Tsa{n}{\bF'} \ar[r]^{\ \simeq} & \bN}$. This yields a quasi-isomorphism
\begin{equation}\label{eq:qi}
\xymatrix@C=0.5cm{\tp{\bT}{\Tsa{n}{\bF'}} \ar[r]^{\ \simeq} & \tp{\bT}{\bN}}
\end{equation}
by \cite[Proposition 2.14]{CFH}. Consider the exact sequence
$$0 \to \tp{\bT}{\Tha{n-1}{\bF'}} \to \tp{\bT}{\Tsa{n}{\bF'}} \to \susp^{n}(\tp{\bT}{\Co[n]{\bF'}}) \to 0.$$
Since the sequence $0 \to \Tha{n-1}{\bF'} \to \bF' \to \Thb{n}{\bF'} \to 0$ is exact, $\Tha{n-1}{\bF'}$ is a semi-flat complex, and so $\tp{\bT}{\Tha{n-1}{\bF'}}$ is acyclic. Thus one gets a quasi-isomorphism $\xymatrix@C=0.5cm{\tp{\bT}{\Tsa{n}{\bF'}} \ar[r]^{\simeq\ \ \ \ \ \ } & \susp^{n}(\tp{\bT}{\Co[n]{\bF'}})}$. Now the desired isomorphism in the statement holds by (\ref{eq:qi}).
\end{prf*}

The next result is proved similarly.

\begin{lem}\label{lem:syzygy dual}
Let $\bT'$ be an acyclic complex of flat $R$-modules, and let $\bM$ be a bounded above $\Rop$-complex with $\bF \qis \bM$ a semi-flat replacement. Then for each $m\geq\sup \bM$ there is an isomorphism $\tp{\bM}{\bT'}\simeq \susp^{m}(\tp{\Co[m]\bF}{\bT'})$ in $\D[\kk]$.
\end{lem}

\begin{lem}\label{thm:tate}
Let $\bM$ be a bounded above $\Rop$-complex with $\bF \qis \bM$ a semi-flat replacement and $\bT$ an acyclic complexes of flat $\Rop$-modules such that for
some $g\in\mathbb{Z}$ there is an isomorphism $\Thb{g}{\bT} \is \Thb{g}{\bF}$. Let $\bN$ be a bounded above $R$-complex with $\bF' \qis \bN$ a semi-flat replacement and $\bT'$ an acyclic complexes of flat $R$-modules such that for
some $g'\in\mathbb{Z}$ there is an isomorphism $\Thb{g'}{\bT'} \is \Thb{g'}{\bF'}$. Then there is an isomorphism $\tp{\bT}{\bN}\simeq\tp{\bM}{\bT'}$ in $\D[\kk]$.
\end{lem}
\begin{prf*}
Set $m=\max\{\sup \bM, g\}$ and $n=\max\{\sup \bN, g'\}$. One gets
\begin{align*}
    \tp{\bT}{\bN}
    &\simeq\susp^{n}(\tp{\bT}{\Co[n]{\bF'}})\\
    &\cong\susp^{n}(\tp{\bT}{\Co[n]{\bT'}})\\
    &\simeq\susp^{m}(\tp{\Co[m]{\bT}}{\bT'})\\
    &\cong\susp^{m}(\tp{\Co[m]{\bF}}{\bT'})\\
    &\simeq\tp{\bM}{\bT'},
\end{align*}
where the first and the last isomorphisms hold by Lemmas \ref{lem:syzygy} and \ref{lem:syzygy dual}, respectively, and the middle one follows from \cite[Lemma 4.1]{CCLP1}.
\end{prf*}

The following balanced result for Tate homology is immediately by Lemma \ref{thm:tate}, where part (b) was first proved by Christensen and Jorgensen in \cite[Theorem 3.7]{CJ} using pinched complexes.

\begin{thm}\label{cor:fb}
Let $\bM$ be a bounded above $\Rop$-complex of finite Gorenstein projective dimension and let $\bN$ be a bounded above $R$-complex.
  \begin{prt}
  \item I\fsp\ $\bN$ has finite Gorenstein flat dimension with $(\bT, \bF)$ a Tate flat resolution, then for each $i\in\ZZ$ there is an isomorphism
$$\Ttor{i}{\bM}{\bN}\cong\H[i]{\tp{\bM}{\bT}}.$$
  \item I\fsp\ $\bN$ has finite Gorenstein projective dimension with $\bS \to \bP \to \bN$ a complete projective resolution, then for each $i\in\ZZ$ there is an isomorphism
$$\Ttor{i}{\bM}{\bN}\cong\H[i]{\tp{\bM}{\bS}}.$$
  \end{prt}
\end{thm}

\section{Stable and unbounded homology of complexes}\label{4}
\noindent
In this section, we give some developments related to Auslander's use of the transpose in terms of unbounded and stable homology modules, and prove Theorem C advertised in the introduction. We first recall several definitions and terminology from \cite{CCLP1, G}.

For an $\Rop$-complex $\bM$ and an $R$-complex $\bN$, the \emph{unbounded tensor product} $\btp{\bM}{\bN}$ is the $\kk$-complex with the degree-$n$ term
$$(\btp{\bM}{\bN})_{n}=\prod_{i\in\ZZ}(M_{i}\otimes_{R}N_{n-i})$$
and the differential defined as in \ref{tensor}. It contains the tensor product $\tp{\bM}{\bN}$ as a subcomplex. The quotient complex $(\btp{\bM}{\bN})/(\tp{\bM}{\bN})$ is called the \emph{stable tensor product}, and it is denoted $\ttp{\bM}{\bN}$.

\begin{dfn}
  \label{dfn:tTor}
  Let $\bM$ be an $\Rop$-complex and $\bN$ an $R$-complex. Let $\bP \qra \bM$
  be a semi-projective resolution and let $\bN \qra \bI$ be a semi-injective
  resolution. For each $i\in\ZZ$, the $i$th \emph{unbounded homology} of $\bM$ and $\bN$ over $R$ is
    $$\bTor{i}{\bM}{\bN} \deq \H[i]{\btp{\bP}{\bI}},$$
  and the $i$th \emph{stable homology} of $\bM$ and $\bN$ over $R$ is
    $$\Stor[R]{i}{\bM}{\bN} \deq \H[i+1]{\ttp{\bP}{\bI}}.$$
\end{dfn}

The proof of the next result is similar to that of \cite[Proposition 2.6]{CCLP1}.

\begin{prp}\label{prp:stable homology}
Let $\bM$ be an $\Rop$-complex with $\bF\qis \bM$ a semi-flat replacement, and let $\bN$ be a homologically bounded above $R$-complex with $\bN\qra \bI$ a semi-injective resolution such that $\sup \bI<\infty$. Then for each $i\in\ZZ$, we have
  \begin{equation*}
    \bTor{i}{\bM}{\bN} \is \H[i]{\btp{\bF}{\bI}} \qand
    \Stor{i}{\bM}{\bN} \is \H[i+1]{\ttp{\bF}{\bI}}\:.
  \end{equation*}
\end{prp}

\begin{lem}\label{lem:btor}
  Let $\bT$ be an F-totally acyclic complex of flat $\Rop$-modules and $\bN$ a bounded above complex of $R$-modules. Then for all integers $i$ and $n$ we have $$\bTor{i}{\Co[n]{\bT}}{\bN}\is\H[i+n-1]{\tp{\Tha{n-1}{\bT}}{\bN}}.$$
  In particular, if $N$ is an $R$-module, then $\bTor{i}{\Co[n]{\bT}}{N}=0$ for each $i\geq 1$.
\end{lem}
\begin{prf*}
  Consider the degree-wise split exact sequence
  \begin{equation*}
    \label{eq:ttt}
    0 \to \Tha{n-1}{\bT} \to \bT \to \Thb{n}{\bT} \to 0\;.
  \end{equation*}
  Let $\bN \qra \bI$ be a semi-injective resolution such that $\sup \bI<\infty$. The complex $\btp{\bT}{\bI}$
  is acyclic by \cite[Proposition 1.7]{CCLP1}, whence there is an isomorphism
  \begin{equation}
    \label{eq:asdf}
    \H{\btp{\Thb{n}\bT}{\bI}} \dis \susp\H{\btp{\Tha{n-1}\bT}{\bI}}\;.
  \end{equation}
  Since the canonical map $\susp^{-n}\Thb{n}{\bT} \to \Co[n]{\bT}$ is a
  flat resolution, one has
  \begin{equation*}
    \begin{aligned}
      \bTor{i}{\Co[n]{\bT}}{\bN} 
      &\is \H[i]{\btp{(\susp^{-n}\Thb{n}{\bT})}{\bI}}\\
      &\is \H[i+n]{\btp{\Thb{n}{\bT}}{\bI}}\\
      &\is \H[i+n-1]{\btp{\Tha{n-1}{\bT}}{\bI}}\\
      &=   \H[i+n-1]{\tp{\Tha{n-1}{\bT}}{\bI}}\\
      &\is \H[i+n-1]{\tp{\Tha{n-1}{\bT}}{\bN}},
    \end{aligned}
  \end{equation*}
  where the first isomorphism holds by Proposition \ref{prp:stable homology}, the third one follows from the isomorphism \eqref{asdf}, and the last one holds by \cite[Proposition 2.14]{CFH}.
\end{prf*}

\begin{lem}\label{prp:Stor-Tateflat}
  Let $\bM$ be an $\Rop$-complex of finite Gorenstein flat dimension with $(\bT,\bF)$ a Tate flat resolution, and let $\bN$ be a bounded above $R$-complex. For every $n\in\ZZ$ there is an
  exact sequence
  $$\xymatrix@C=0.22cm{
    \cdots \ar[r] & \Stor{n+1}{\bM}{\bN} \ar[r] & \Tor{1}{\Co[n]{\bT}}{\bN} \ar[r] & \bTor{1}{\Co[n]{\bT}}{\bN} \ar[r] & \Stor{n}{\bM}{\bN} \ar[r] & \cdots}.$$
\end{lem}
\begin{prf*}
  Consider the degree-wise split exact sequence
  \begin{equation}
    \label{eq:ttt}
    0 \to \Tha{n-1}{\bT} \to \bT \to \Thb{n}{\bT} \to 0\;.
  \end{equation}
  Then by Lemma \ref{lem:btor}, one has
  \begin{equation}
    \label{eq:btorb1}
    \H[i]{\tp{\Tha{n-1}{\bT}}{\bN}} \is \bTor{i-n+1}{\Co[n]{\bT}}{\bN}
  \end{equation}
  for each $i\in\ZZ$. It is clear that
  \begin{equation}
    \label{eq:btorb}
    \H[i]{\tp{\Thb{n}{\bT}}{\bN}} \is \Tor{i-n}{\Co[n]{\bT}}{\bN}\;.
  \end{equation}
  From \eqref{ttt} one gets the exact sequence
  \begin{equation*}
    0 \to \tp{\Tha{n-1}{\bT}}{\bN} \to \tp{\bT}{\bN} \to \tp{\Thb{n}{\bT}}{\bN} \to 0\;,
  \end{equation*}
  and in view of \cite[Theorem 3.10]{CCLP1} and the isomorphisms
  \eqref{btorb1} and \eqref{btorb} the associated sequence in homology is
   $$\xymatrix@C=0.22cm{
    \cdots \ar[r] & \Stor{n+1}{\bM}{\bN} \ar[r] & \Tor{1}{\Co[n]{\bT}}{\bN} \ar[r] & \bTor{1}{\Co[n]{\bT}}{\bN} \ar[r] & \Stor{n}{\bM}{\bN} \ar[r] & \cdots}$$
  as desired.
\end{prf*}

\begin{thm}
  \label{cor:sequence}
  Let $\bM$ be an $\Rop$-complex of finite Gorenstein flat dimension with $(\bT,\bF)$ a Tate flat resolution, and let $N$ be an $R$-module. For every $n\in\ZZ$ there is an
  exact sequence
  \begin{equation*}
    0 \to \Stor{n}{\bM}{N} \to \tp{\Co[n]{\bT}}{N} \to \bTor{0}{\Co[n]{\bT}}{N}
    \to \Stor{n-1}{\bM}{N} \to 0\;.
  \end{equation*}
  In particular, if $M$ is a Gorenstein flat $\Rop$-module, then there is an
  exact sequence
    \begin{equation*}
    0 \to \Stor{0}{M}{N} \to \tp{M}{N} \to \bTor{0}{M}{N}
    \to \Stor{-1}{M}{N} \to 0\;.
  \end{equation*}
\end{thm}
\begin{prf*}
We notice that $\bTor{1}{\Co[n]{\bT}}{N}=0$ by Lemma \ref{lem:btor}. The desired sequence in the statement now follows from Lemma \ref{prp:Stor-Tateflat}.
\end{prf*}

\begin{prp}\label{prp:gproj}
  Let $R$ be a noetherian ring, and let $M$ be a finitely generated Gorenstein projective $\Rop$-module and $\bN$ a bounded above complex of $R$-modules. Then for each $i\in\ZZ$, there is an isomorphism
  $$\bTor{i}{M}{\bN}\is\Ext{-i}{\Hom[\Rop]{M}{R}}{\bN}.$$
  In particular, if $N$ is an $R$-module, then $\bTor{0}{M}{N}\is\Hom{\Hom[\Rop]{M}{R}}{N}$ and $\bTor{i}{M}{N}=0$ for each $i\geq 1$.
\end{prp}
\begin{prf*}
  There exists an acyclic complex $\bT$ of finitely generated free $\Rop$-modules, such that $\Hom[\Rop]{\bT}{R}$ is acyclic and $\Co[0]{\bT}\is M$. Note that $\Hom[\Rop]{\bT}{R}$ is a complex of finitely generated free $\Rop$-modules. Then by \cite[Theorem 4.5.7]{CFoH}, one gets that
    $\tp{\bT}{E} \is
    \tp{\Hom{\Hom[\Rop]{\bT}{R}}{R}}{E} \is \Hom{\Hom[\Rop]{\bT}{R}}{E}$
  is acyclic for each injective $R$-module $E$, and so $\bT$ is an F-totally acyclic complex. Thus for every $i\in\ZZ$ one has
  \begin{align*}
    \bTor{i}{M}{\bN}
    &\dis \H[i-1]{\tp{\Tha{-1}{\bT}}{\bN}}\\
    &\dis \H[i-1]{\tp{\Hom{\Hom[\Rop]{\Tha{-1}{\bT}}{R}}{R}}{\bN}}\\
    &\dis \H[i-1]{\Hom{\Hom[\Rop]{\Tha{-1}{\bT}}{R}}{\bN}}\\
    &\dis \H[i]{\Hom{\susp^{-1}\Hom[\Rop]{\Tha{-1}{\bT}}{R}}{\bN}}\\
    &\dis \Ext{-i}{\Hom[\Rop]{M}{R}}{\bN}\;,
  \end{align*}
  where the first isomorphism holds by Lemma \ref{lem:btor}, the third one follows from \cite[Theorem 4.5.7]{CFoH} as $\bN$ is bounded above, and the last one holds since
  $$\susp^{-1}\Hom[\Rop]{\Tha{-1}{\bT}}{R} \to \Hom[\Rop]{M}{R}$$
  is a projective resolution.
\end{prf*}

If $R$ is a noetherian ring and $M$ is a finitely generated Gorenstein projective $\Rop$-module, then the stable homology $\Stor{*}{M}{-}$ coincides with the Tate homology $\Ttor{*}{M}{-}$ by \cite[Theorem 6.4]{CCLP1}. Thus the following result follows from Theorem \ref{cor:sequence} and Proposition \ref{prp:gproj}, which was proved by Christensen and Jorgensen over a commutative noetherian local ring; see \cite[Proposition 6.3]{CJ}.

\begin{cor}\label{CJ}
  Let $R$ be a noetherian ring, and let $M$ be a finitely generated Gorenstein projective $\Rop$-module and $N$ an $R$-module. Then there is an exact sequence
    \begin{equation*}
    0 \to \Ttor{0}{M}{N} \to \tp{M}{N} \to \Hom{\Hom[\Rop]{M}{R}}{N}
    \to \Ttor{-1}{M}{N} \to 0\;.
  \end{equation*}
\end{cor}

\section{Unbounded homology of complexes of finite Gorenstein dimension}\label{5}
\noindent
In this section, we further study unbounded homology for complexes of finite Gorenstein flat dimension, and we compare unbounded homology to relative homology in the case of module arguments.

Let $R$ be a left coherent ring, and let $\bM$ be an $\Rop$-complex of finite Gorenstein flat dimension. Recall from Theorem \ref{prp:resolution} that $\bM$ admits a complete flat resolution $\tau: \bT\to \bF$ such that $\tau_i$ is split surjective. Set $\bK=\Ker\tau$. Then we have the following result.

\begin{thm}\label{prp:unbounded homology}
Let $\bM$ be an $\Rop$-complex admitting a complete flat resolution $\tau: \bT\to \bF$ such that $\tau_i$ is split surjective, and let $\bN$ be a homologically bounded above $R$-complex with $\bN\qra \bI$ a semi-injective resolution such that $\sup \bI<\infty$. Then for each $i\in\ZZ$, there are isomorphisms
  \begin{equation*}
    \bTor{i}{\bM}{\bN}\is\H[i-1]{\tp{\Ker{\tau}}{\bI}} \qand
    \Stor{i}{\bM}{\bN}\is\H[i]{\tp{\bT}{\bI}}\:.
  \end{equation*}
If furthermore $\bN$ is a bounded above complex, then there are isomorphisms
  \begin{equation*}
    \bTor{i}{\bM}{\bN}\is\H[i-1]{\tp{\Ker{\tau}}{\bN}} \qand
    \Stor{i}{\bM}{\bN}\is\H[i]{\tp{\bT}{\bN}}\:.
  \end{equation*}
\end{thm}
\begin{prf*}
Set $\bK=\Ker{\tau}$. Consider the degree-wise split exact sequence $0\to \bK \to \bT \to \bF \to 0$. By \cite[1.5(c)]{CCLP1} one has the following commutative diagram with exact rows and columns:
$$\xymatrix@C=20pt@R=20pt{
  & 0 \ar[d] & 0 \ar[d] & 0 \ar[d] & \\
  0 \ar[r] & \tp{\bK}{\bI} \ar[d]\ar[r] & \tp{\bT}{\bI} \ar[d]\ar[r] & \tp{\bF}{\bI} \ar[d]\ar[r] &  0 \\
  0 \ar[r] & \btp{\bK}{\bI} \ar[d]\ar[r] & \btp{\bT}{\bI} \ar[d]\ar[r] & \btp{\bF}{\bI} \ar[d]\ar[r] &  0 \\
  0 \ar[r] & \ttp{\bK}{\bI} \ar[d]\ar[r] & \ttp{\bT}{\bI} \ar[d]\ar[r] & \ttp{\bF}{\bI} \ar[d]\ar[r] &  0\ . \\
  & 0 & 0 & 0  & }$$

Since $\sup \bI<\infty$, $\btp{\bT}{\bI}$ is acyclic by \cite[Proposition 1.7(a)]{CCLP1}. We observe that $\bK$ is bounded above. Then $\ttp{\bK}{\bI}=0$, and so by Proposition \ref{prp:stable homology}, for each $i\in\ZZ$, one has
$$\bTor{i}{\bM}{\bN}\is\H[i]{\btp{\bF}{\bI}}\is\H[i-1]{\btp{\bK}{\bI}}=\H[i-1]{\tp{\bK}{\bI}},$$
and
$$\Stor{i}{\bM}{\bN}\is\H[i+1]{\ttp{\bF}{\bI}}=\H[i+1]{\ttp{\bT}{\bI}}\is\H[i]{\tp{\bT}{\bI}}.$$

If furthermore $\bN$ is a bounded above complex, then there are quasi-isomorphisms $\tp{\bK}{\bN}\qra\tp{\bK}{\bI}$ and $\tp{\bT}{\bN}\qra\tp{\bT}{\bI}$ by \cite[Proposition 2.14]{CFH}, as both $\bK$ and $\bT$ are complexes of flat $\Rop$-modules. Thus for each $i\in\ZZ$, one has $\bTor{i}{\bM}{\bN}\is\H[i-1]{\tp{\bK}{\bN}}$ and $\Stor{i}{\bM}{\bN}\is\H[i]{\tp{\bT}{\bN}}$.
\end{prf*}

\begin{cor}\label{cor:vanishing}
Let $R$ be a left coherent ring, and let $\bM$ be an $\Rop$-complex of finite Gorenstein flat dimension and $\bN$ a homologically bounded above $R$-complex. Then $\bTor{i}{\bM}{\bN}=0$ for each $i>\Gfd[\Rop]\bM+\sup\H{\bN}$.
\end{cor}
\begin{prf*}
Let $\Gfd[\Rop]\bM=g<\infty$. By Theorem \ref{prp:resolution} $\bM$ admits a complete flat resolution $\tau: \bT\to \bF$ such that $\tau_i$ is split surjective for each $i\in\ZZ$ and $\tau_i$ is an isomorphism for each $i\geq g$. Set $\bK=\Ker\tau$. Then one has $\sup \bK=g-1$. Fix a semi-injective resolution $\bN\qra \bI$ such that $\sup \bI=\sup\H{\bN}$. Then $\sup(\tp{\bK}{\bI})\leq\sup \bK+\sup \bI=g-1+\sup\H{\bN}$. Thus by Theorem \ref{prp:unbounded homology}, $\bTor{i}{\bM}{\bN}\is\H[i-1]{\tp{\bK}{\bI}}=0$ for each $i>g+\sup\H{\bN}$.
\end{prf*}

The next result compares to Proposition \ref{prp:gfd}.

\begin{cor}\label{cor:gflatdim}
Let $R$ be a left coherent ring, and let $\bM$ be an $\Rop$-complex of finite Gorenstein flat dimension. Then one has
  \begin{align*}
    \Gfd[\Rop]{\bM}
    & \deq \sup\left\{ n \in \ZZ \:
      \left|
        \begin{array}{c}
          \bTor{n+\sup\H{\bN}}{\bM}{\bN}\neq 0\\
          \text{for some $R$-complex $\bN$ with}\\
          \text{$\sup\H{\bN}<\infty$}
        \end{array}
      \right.
    \right\}\\
    & \deq \sup\left\{ n \in \ZZ \:
      \left|
        \begin{array}{c}
          \bTor{n}{\bM}{N}\neq 0\\
          \text{for some $R$-module $N$}
        \end{array}
      \right.
    \right\}.
  \end{align*}
\end{cor}
\begin{prf*}
Let $\Gfd[\Rop]{\bM}=g<\infty$. If $g=-\infty$, then $\bM$ is acyclic; see Remark \ref{rmk:gfd-acyclic}. If fix a semi-projective resolution $\bP\qra \bM$, then $\bP$ is contractible, and so $\btp{\bP}{\bI}$ is acyclic for each $R$-complex $\bI$ by \cite[1.5]{CCLP1}. Thus $\bTor{n}{\bM}{-}=0$ for each $n\in\ZZ$, and so the equalities in the statement hold. We let $g\in \ZZ$, and let $s$ (resp., $s'$) denote the number on the right side of the first (resp., second) equality in the statement. Corollary \ref{cor:vanishing} yields $g\geq s$. Obviously, $s\geq s'$, so it remains to prove $s'\geq g$. By Proposition \ref{prp:gfd}, there is an injective $R$-module $E$ such that $\Tor{g}{\bM}{E}\neq 0$. Since $\Stor{i}{\bM}{E}=0$ for each $i\in\ZZ$, one has $\bTor{g}{\bM}{E}\is\Tor{g}{\bM}{E}\neq 0$; see \cite[(2.5)]{CCLP1}. This yields $s'\geq g$.
\end{prf*}

The following is a balanced result for unbounded homology of complexes.

\begin{thm}\label{thm:balance}
Let $R$ be a coherent ring (that is, both left and right coherent), and let $\bM$ be an $\Rop$-complex and $\bN$ an $R$-complex, both of which are of finite Gorenstein flat dimension. Then for each $i\in\ZZ$ there is an isomorphism
$$\bTor{i}{\bM}{\bN}\is\bTor[\Rop]{i}{\bN}{\bM}.$$
\end{thm}
\begin{prf*}
Let $\Gfd[\Rop]\bM=s$ and $\Gfd \bN=t$. Then by \cite[Proposition 5.7]{CKL} one has $\sup\H{\bM}\leq s<\infty$ and $\sup\H{\bN}\leq t<\infty$. It follows from Theorem \ref{prp:resolution} that $\bM$ admits a complete flat resolution $\tau: \bT\to \bF$ such that $\tau_i$ is split surjective for each $i\in\ZZ$ and $\tau_i$ is an isomorphism for each $i\geq s$, and $\bN$ admits a complete flat resolution $\tau': \bT'\to \bF'$ such that $\tau_i'$ is split surjective for each $i\in\ZZ$ and $\tau_i'$ is an isomorphism for each $i\geq t$. Fix semi-injective resolutions $\bM\qra \bI$ and $\bN\qra \bI'$ such that $\sup \bI=\sup\H{\bM}$ and $\sup \bI'=\sup\H{\bN}$. Set $\bK=\Ker\tau$ and $\bK'=\Ker\tau'$. Then for each $i\in\ZZ$, one has
  \begin{align*}
    \bTor{i}{\bM}{\bN}
    &\is \H[i-1]{\tp{\bK}{\bI'}}\\
    &=   \H[i-1]{\tp{\Thb{i-2-t}{\bK}}{\bI'}}\\
    &\is \H[i-1]{\tp{\Thb{i-2-t}{\bK}}{\bF'}}\\
    &\is \H[i-2]{\tp{\Thb{i-2-t}{\bK}}{\bK'}}\\
    &=   \H[i-2]{\tp{\bK}{\bK'}}\\
    &\is \H[i-2]{\tp[\Rop]{\bK'}{\bK}}\\
    &=   \H[i-2]{\tp[\Rop]{\Thb{i-2-s}{\bK'}}{\bK}}\\
    &\is \H[i-1]{\tp[\Rop]{\Thb{i-2-s}{\bK'}}{\bF}}\\
    &\is \H[i-1]{\tp[\Rop]{\Thb{i-2-s}{\bK'}}{\bI}}\\
    &=   \H[i-1]{\tp[\Rop]{\bK'}{\bI}}\\
    &\is \bTor[\Rop]{i}{\bN}{\bM}\:,
  \end{align*}
  where the first and the last isomorphisms follow from Theorem \ref{prp:unbounded homology}, the four equations hold by Lemma \ref{lem:bounded} since $\sup \bI'\leq t$, $\sup \bK'=t-1$, $\sup \bK=s-1$ and $\sup \bI\leq s$, the remaining isomorphisms hold since both $\Thb{i-2-t}{\bK}$ and $\Thb{i-2-s}{\bK'}$ are semi-flat complexes.
\end{prf*}

\begin{ipg}
Let $M$ be an $\Rop$-module of finite Gorenstein projective dimension. Then there exists a proper Gorenstein projective resolution of $M$, that is, a quasi-isomorphism $\pi: \bG\qra M$ with $G_{i}$ Gorenstein projective for each $i\geq 0$ and $G_{i}=0$ for each $i<0$, such that $\Hom[\Rop]{L}{\Cone\pi}$ is acyclic for each Gorenstein projective $\Rop$-module $L$; see \cite[Theorem 2.10]{H1}. Following Holm \cite[Section 4]{H2}, for each $R$-module $N$ and each $i\geq 0$, the $i$th \emph{relative homology} of $M$ and $N$ over $R$ is defined as  $$\mathrm{Tor}^{\mathcal{GP}}_{i}(M,N)=\H[i]{\tp{\bG}{N}}.$$
If $R$ is a left coherent ring and $M$ is an $\Rop$-module of finite Gorenstein flat dimension, then by \cite[Theorem 3.23]{H1}, $M$ admits a proper Gorenstein flat resolution $\bG'\qra M$. Following \cite[Section 4]{H2}, for each $R$-module $N$ and each $i\geq 0$, the $i$th \emph{relative homology based on flats} of $M$ and $N$ over $R$ is defined as  $$\mathrm{Tor}^{\mathcal{GF}}_{i}(M,N)=\H[i]{\tp{\bG'}{N}}.$$
\end{ipg}

\begin{ipg}
  We recall the invariant
    $\splfR \deq \sup\setof{\pd{F}}{F\ \text{\rm is a flat $R$-module}}$. Since an arbitrary direct sum of flat $R$-modules is flat, the
  invariant $\splfR$ is finite if and only if every flat $R$-module
  has finite projective dimension.
  If $R$ is commutative noetherian of finite Krull dimension $d$, then one has
  $\splfR \le d$ by Jensen \cite[Proposition 6]{J}. Osofsky
  \cite[3.1]{O} gives examples of rings for which the splf
  invariant is infinite.
\end{ipg}

We end this section with the following results that establish the relation between unbounded homology groups and relative homology groups.

\begin{prp}\label{prp:grelative}
Let $R$ be a left coherent ring with $\mathrm{splf}\Rop<\infty$, and let $M$ be an $\Rop$-module of finite Gorenstein projective dimension and $N$ an $R$-module. Then there is an isomorphism
$$\bTor{i}{M}{N}\is\mathrm{Tor}^{\mathcal{GP}}_{i}(M,N)$$
for each $i\geq 2$, and an exact sequence
$$0 \to \mathrm{Tor}^{\mathcal{GP}}_{1}(M,N) \to \bTor{1}{M}{N} \to \Stor{0}{M}{N}.$$
\end{prp}
\begin{prf*}
By \cite[Theorem 3.4]{V}, $M$ admits a complete projective resolution $\tau: \bT\to \bP$ such that $\tau_i$ is split surjective for each $i\in\ZZ$. This yields a degree-wise split exact sequence
$$0 \to \susp^{-1}\bG \to \widetilde{\bT} \to \bP \to 0$$
 such that $\bG\qra M$ is a proper Gorenstein projective resolution and $\Thb{0}{(\susp^{-1}\bG)}=\Thb{0}{(\Ker\tau)}$, and $\widetilde{\bT}$ is acyclic with $\Thb{0}{\widetilde{\bT}}=\Thb{0}{\bT}$ and $\widetilde{\bT}_{i}=0$ for each $i\leq -2$; see \cite[(3.8)]{AM}. Since all Gorenstein projective $\Rop$-modules are Gorenstein flat by the proof of \cite[Proposition 3.4]{H1}, $\tau: \bT\to \bP$ is a complete flat resolution of $M$, and so by Theorem \ref{prp:unbounded homology}, for each $i\geq 2$, one has
  \begin{align*}
    \bTor{i}{M}{N}
    &\is \H[i-1]{\tp{\Ker\tau}{N}}\\
    &=   \H[i-1]{\tp{\susp^{-1}\bG}{N}}\\
    &\is \H[i]{\tp{\bG}{N}}\\
    &\is \mathrm{Tor}^{\mathcal{GP}}_{i}(M,N)\:.
  \end{align*}
  Consider the exact sequence
  $$0 \to \tp{\susp^{-1}\bG}{N} \to \tp{\widetilde{\bT}}{N} \to \tp{\bP}{N} \to 0.$$
  Then one gets the exact sequence
  $$\H[1]{\tp{\widetilde{\bT}}{N}} \to \H[1]{\tp{\bP}{N}} \to \H[0]{\tp{\susp^{-1}\bG}{N}} \to \H[0]{\tp{\widetilde{\bT}}{N}}.$$
  It is clear that $\H[0]{\tp{\widetilde{\bT}}{N}}=0$ and $\H[1]{\tp{\widetilde{\bT}}{N}}=\H[1]{\tp{\bT}{N}}$, and
  $$\H[0]{\tp{\susp^{-1}\bG}{N}} \is \mathrm{Tor}^{\mathcal{GP}}_{1}(M,N).$$
  Thus one has
  \begin{equation}\label{1}
  \mathrm{Tor}^{\mathcal{GP}}_{1}(M,N) \is \Coker{\H[1]{\tp{\tau}{N}}}.
  \end{equation}
  On the other hand, consider the exact sequence
  $$\xymatrix@C=0.5cm{0 \ar[r] & \tp{(\Ker\tau)}{N} \ar[r] & \tp{\bT}{N} \ar[rr]^{\tp{\tau}{N}} && \tp{\bP}{N} \ar[r] & 0.}$$
  Then one gets the exact sequence
  \begin{equation}\label{2}
  \xymatrix@C=0.5cm{\H[1]{\tp{\bT}{N}} \ar[rr]^{\H[1]{\tp{\tau}{N}}} && \H[1]{\tp{\bP}{N}} \ar[r] & \H[0]{\tp{(\Ker\tau)}{N}} \ar[r] & \H[0]{\tp{\bT}{N}}.}
  \end{equation}
  Since $\H[0]{\tp{(\Ker\tau)}{N}} \is \bTor{1}{M}{N}$ and $\H[0]{\tp{\bT}{N}} \is \Stor{0}{M}{N}$ by Theorem \ref{prp:unbounded homology} and \cite[Theorem 3.10]{CCLP1}, (\ref{1}) and (\ref{2}) yield the exact sequence
  \begin{equation*}
    0 \to \mathrm{Tor}^{\mathcal{GP}}_{1}(M,N) \to \bTor{1}{M}{N} \to \Stor{0}{M}{N}\:.\qedhere
  \end{equation*}
\end{prf*}

\begin{prp}\label{prp:grelativeflat}
Let $R$ be a left coherent ring, and let $M$ be a cotorsion $\Rop$-module of finite Gorenstein flat dimension and $N$ an $R$-module. Then there is an isomorphism
$$\bTor{i}{M}{N}\is\mathrm{Tor}^{\mathcal{GF}}_{i}(M,N)$$
for each $i\geq 2$, and an exact sequence
$$0 \to \mathrm{Tor}^{\mathcal{GF}}_{1}(M,N) \to \bTor{1}{M}{N} \to \Stor{0}{M}{N}.$$
\end{prp}
\begin{prf*}
By \cite[Lemma 4.2]{L1}, $M$ admits a complete flat resolution $\tau: \bT\to \bF$ such that for each $i\in\ZZ$, $T_i$ and $F_i$ are cotorsion, and $\tau_i$ is split surjective. This yields a degree-wise split exact sequence
$$0 \to \susp^{-1}\bG \to \widetilde{\bT} \to \bF \to 0$$
 such that $\bG\qra M$ is a proper Gorenstein flat resolution and $\Thb{0}{(\susp^{-1}\bG)}=\Thb{0}{(\Ker\tau)}$, and $\widetilde{\bT}$ is acyclic with $\Thb{0}{\widetilde{\bT}}=\Thb{0}{\bT}$ and $\widetilde{\bT}_{i}=0$ for each $i\leq -2$; see \cite[Lemma 4.3]{L1}. Now one can complete the proof using the same argument as in the proof of Proposition \ref{prp:grelative}.
\end{prf*}

\section*{Acknowledgments}
\noindent
The author thanks the referee for valuable comments that have improved the presentation at several points.

 \providecommand{\arxiv}[2][AC]{\mbox{\href{http://arxiv.org/abs/#2}{\sf
  arXiv:#2 [math.#1]}}} \def\cprime{$'$}
\providecommand{\bysame}{\leavevmode\hbox to3em{\hrulefill}\thinspace}

\end{document}